\documentclass[12pt,reqno]{amsart}
\usepackage[utf8]{inputenc}
\usepackage{color}
\usepackage[shortlabels]{enumitem}
\usepackage{float}
\usepackage{mathtools}
\usepackage{graphicx}
\usepackage{colortbl}
\usepackage{textcomp}
\usepackage{amssymb}
\usepackage{amsmath}
\usepackage{mathrsfs}
\usepackage{amssymb,amsmath,mathrsfs,latexsym,color,wasysym,graphicx}
\usepackage{enumitem}
\usepackage{multicol}

\theoremstyle{definition}
\newtheorem{thm}{Theorem}[section]
\newtheorem{prop}[thm]{Proposition}
\newtheorem{cor}[thm]{Corollary}
\newtheorem{lemma}[thm]{Lemma}

\newtheorem{rema}[thm]{Remark}
\newtheorem{ex}[thm]{Example}
\newtheorem{defn}[thm]{Definition}

\newcommand{\pf}{\noindent \textit{Proof}:\ }

\def\qed{\ensuremath\square}

\begin{document}

\title[Lie Bracket Algebras over Leavitt Path Algebras]{A Talented Monoid View on Lie Bracket Algebras over Leavitt Path Algebras}
%
\author[W.~Bock]{Wolfgang Bock}
\address{Wolfgang Bock: Technomathematics Group,Department of Mathematics, Univer-
	sity of Kaiserslautern, 67653 Kaiserslautern, Germany}
\email{bock@mathematik.uni-kl.de}

\author[A.~Sebandal]{Alfilgen N. Sebandal}
\address{Alfilgen N. Sebandal: MSU-IIT Iligan, Andres Bonifacio Avenue, Tibanga,9200 Iligan City, Philippines}
\email{alfilgen.sebandal@g.msuiit.edu.ph}

\author[J.~Vilela]{Jocelyn P. Vilela}
\address{Jocelyn P. Vilela: MSU-IIT Iligan, Andres Bonifacio Avenue, Tibanga,9200 Iligan City, Philippines}
\email{jocelyn.vilela@g.msuiit.edu.ph}

%
\date{\today}


\maketitle

\begin{abstract}
In this article, we study properties as simplicity, solvability and nilpotency for Lie bracket algebras arising from Leavitt path algebras, based on the talented monoid of the underlying graph. We show that graded simplicity and simplicity of the Leavitt path algebra can be connected via the Lie bracket algebra. Moreover, we use the Gelfand-Kirillov dimension for the Leavitt path algebra for a classification of nilpotency and solvability. 
\end{abstract}
\maketitle

\section{Introduction}
Leavitt path algebras have become an object of intensive research in recent years, combining both a complete description of the algebraic structure via the geometry of the associated graphs and a complete invariant for the classification of the algebras. For a good overview of recent results, we refer to the overview paper \cite{Abramsdecade} and the monograph \cite{LPAbook}.

A \emph{directed graph} $E$ is a tuple $(E^0, E^1, r,s)$ where $E^0$ and $E^1$ are sets and $r,s$ are maps from $E^1$ to $E^0$.  The elements of $E^0$ are called $vertices$ and the
elements of $E^1$ $edges$. A graph $E$ is \emph{finite} if $E^0$ and $E^1$ are both finite. The edge $e\in E^1$ is viewed as an edge pointing from $s(e)$ to $r(e)$.
As a convention, a (finite) \emph{path} $p$ in $E$ is a sequence
$p=\alpha_1\alpha_2\cdots \alpha_n$
of edges $\alpha_i$ in $E$ such that $r(\alpha_i)=s(\alpha_{i+1})$ for $1\leq i \leq i-1$. Throughout this paper, graph shall mean directed graph.~~ 

Let $E$ be a graph. If $s^{-1}(v)$ is a finite set for every $v\in E^0$, the the graph is called  $row$-$finite$. A vertex $v$ for which $s^{-1}(v)=\varnothing$ is called a $sink$, while a vertex $v$ for which $r^{-1}(v)=\varnothing$ is called a $source$. In other words, $v$ is a sink (resp., source) if $v$ is not a source (resp., range) of any edge of $E$. A vertex which is both a source and a sink is called $isolated$. A vertex $v$ such that $|s^{-1}(v)|$ is infinite is called an $infinite$ $emitter$. If $v$ is either a sink or an infinite emitter, we call $v$ a $singular$ $vertex$; otherwise, $v$ is called a $regular$ $vertex$. Let Sink($E$), Source($E$), Reg($E$), and Inf($E$) denote the set of sinks, sources, regular vertices and infinite emitters of $E$, respectively. A subset $H\subseteq E^0$ is said to be $hereditary$ if for any $e\in E^1$, $s(e)\in H$ implies $r(e)\in H$.  A subset $H\subseteq E^0$ is called $saturated$ if whenever for a regular vertex $v$, $r(s^{-1}(v))\subseteq H$ then $v\in H$. 
\\
\indent  Given a graph $E$, the $covering$ $graph$ $\overline{E}$ of $E$ is given by ${\overline{E}}^0=\{v_n: v\in E^0,n\in \mathbb{Z}   \}$,
${\overline{E}}^1=\{e_n: e\in E^1,n\in \mathbb{Z}   \}$, $s(e_n)=s(e)_n$ and $r(e_n)=r(e)_{n+1}$. For hereditary saturated subsets $H_1$ and $H_2$ of $E$ with $H_1\subseteq H_2$, define the $quotient$ $graph$ $H_2/H_1$ as a graph such that $(H_2/H_1)^0=H_2\setminus H_1$ and $(H_2/H_1)^1=\{e\in E^1: s(e)\in H_2, r(e)\not \in H_1  \}$. The source and range maps of $H_2/H_1$ are restricted from the graph $E$. The $double$ $graph$ $of$ $E$ is defined as the new graph $\widehat{E}=(E^0, E^1\cup (E^1)^*, r',s'  )$, where $(E^1)^*= \{ e^*:e\in E^1 \}$, and the functions $r'$ and $s'$ are defined as
$r'|_{E^1}=r$, $s'|_{E^1}=s$, $r'(e^*)=s(e)$, and $s'(e^*)=r(e)$ for all $e\in E^1$. In other words, each $e^*\in (E^1)^*$ has orientation the reverse of that of its counterpart $e\in E^1$. The elements $(E^1)^*$ are called $ghost$ $edges$.\\
\indent For a graph $E$ and a ring $R$ with identity, we define the \textit{Leavitt path algebra} of $E$, denoted by $L_R(E)$, to be the algebra generated by the sets $\{v:v\in E^0\}$, $\{ \alpha :\alpha \in E^1  \}$ and $\{ \alpha^* : \alpha \in E^1  \}$ with coefficients in $R$, subject to the relations
\begin{enumerate}
\item[\textnormal{(i)}]
$v_iv_j=\delta_{i,j}v_i$ for every $v_i, v_j\in E^0$;
\item[\textnormal{(ii)}] $s(\alpha)\alpha=\alpha = \alpha r(\alpha)$ and $r(\alpha)\alpha^*=\alpha^*=\alpha^*s(\alpha)$ for all $\alpha \in E^1$;

\item[\textnormal{(iii)}] $\alpha^*\alpha'=\delta_{\alpha,\alpha'}r(\alpha)$ for all $\alpha, \alpha'\in E^1$;

\item[\textnormal{(iv)}]$\sum_{ \{\alpha \in E^1 :s(\alpha)=v    \}  } \alpha \alpha^*=v$ for every $v\in \textnormal{Reg}(E)$.

\end{enumerate}
This article focuses on Leavitt path algebras with coefficients in a field $K$. Already in the first paper in 2005, a simplicity criterion  for this algebra was given \cite{AP05}. In a second paper in 2007, the graph monoid of the associated graph was linked to $K$-theoretic properties of the algebras \cite{AMP07}. 

Using the vertices of a graph, one can create a free monoid with additional relations, called a graph monoid. The  graded version of this graph monoid, the so called talented monoid $T_E$, was introduced in this form by Hazrat and Li in 2020 \cite{A1}. Based on the talented monoid, several algebraic properties of the corresponding Leavitt path algebra, such as simplicity, purely infinite simplicity, or the lattice of ideal could be completely retrieved. In \cite{hazratGK}, the Gelfand-Kirillov dimension for the special case, namely graphs with disjoint cycles, has been classified via the talented monoid. 

A \emph{Lie algebra} over a field $K$ is a $K$-vector space $L$ together with a mapping 
$[-,-]:L\times L\longrightarrow L$
such that
\begin{enumerate}
    \item [(i)] $[-,-]$ is bilinear,
\item [(ii)] $[x,x]=0$ for all $x\in L$,
\item [(iii)] $[x,[y,z]]+[y,[z,x]]+[[z,[x,y]]]=0$ for all $x,y,z\in L$.
\end{enumerate} For an associative $K$-algebra $R$, one can construct the corresponding Lie bracket algebra $[R,R]$ of $R$, defined as $\mathrm{span}\{xy-yx: x,y \in R\}$. Note that the resulting Lie algebra may be not associative. Similar constructions can be done if $R$ is a Novikov algebra, see \cite{Shestakov}. In particular, this is possible for Leavitt path algebras. The Lie bracket algebras arising from Leavitt path algebras have been investigated by \cite{abramsmesyan,Alahmedi} for simplicity and \cite{namzerui} for solvability and nilpotency. This paper shall be a translation of these concepts into the talented monoid point of view.

In this article, we focus on two properties of Lie algebras namely: simplicity conditions of the Lie algebra and solvability or nilpotency conditions.
After a brief overview on preliminaries and previous results, the article's main results are based on, in Section 4, we translate the results of \cite{namzerui} for solvabilty and nilpotency in the context on the talented monoid. Indeed we show that the Gelfand-Kirillov dimension of the Leavitt path algebra over a finite graph has to be bounded by one for the Leavitt path algebra to be Lie solvable. In Theorem \ref{2thm7}, we give - based on the characteristic of $K$, a full description for Lie solvability based on the Gelfand-Kirillov dimension  and ideals of the talented monoid. 
Section 5 is dedicated to the simplicity of the Lie bracket algebra over a Leavitt path algebra. For this we use results from \cite{abramsmesyan} for the case where $L_K(E)$ is simple and \cite{Alahmedi}, in the non-simple case. After a classification of a balloon in the language of a talented monoid, in Theorem \ref{4thm3}, we deduce the relationship of Lie simplicity to cyclic composition series on $T_E$. Based on the relations on Lie simplicity to saturated hereditary sets, we deduce that Lie simplicity and graded simplicity of $L_K(E)$ are equivalent to $L_K(E)$ being simple and $1_{L_{K}(E)}\notin [L_K(E), L_K(E)]$.\\
The latter could lead to new ideas to prove the graded isomorphism conjecture \cite{A1_8, Ara} and also \cite{Vas, Vas2} by a detour via Lie bracket algebras.

\section{Preliminaries}

\begin{defn} \cite{monoid.q}
	Let $M$, $M_1$ and $M_2$ be commutative monoids. 
	\begin{enumerate}
		\item[\textnormal{(i)}] For any submonoid $H$ of $M$,  we define a binary relation $\rho_H$ in $M$ by 
		$x\rho_Hy$ if~and~only~if $(x+H)\cap (y+H)\neq \varnothing.$
		\item[\textnormal{(ii)}]  For a mapping $f:M_1\rightarrow M_2$, define a relation 
		$x\rho_f y$ if and only if $f(x)=f(y).$
	\end{enumerate}
\end{defn}
It can be shown that both $\rho_f$ and $\rho_H$ are equivalence relations for any mapping $f$ and submonoid $H$. For any submonoid $H$ of $M$, the set
$${M}/{H} \cong {M}/{\rho_H} =\lbrace \rho_{H}(x):x\in M\rbrace \hfill$$

On a set $X$ with an equipped preordering $\preccurlyeq$, two elements $x$ and $y$ are said to be \emph{comparable} if $x\preccurlyeq y$ or $y\preccurlyeq x$. There is a natural algebraic pre-ordering on a commutative monoid $M$ defined by $a\leq b$ if $b=a+c$, for some $c\in M$. Throughout, $a\parallel b$ shall mean the elements $a$ and $b$ are not comparable.

\begin{defn} \cite{A1}
 Let $M$ be a monoid with a group $\Gamma$ acting on it. Then $M$ is said to be a  $\Gamma$-\emph{monoid} and the action on an element $a\in M$ by $\alpha \in \Gamma$ shall be denoted by $^\alpha a$.
	Let $M, M_1$ and $M_2$ be monoids. $\Gamma$ a group acting on $M, M_1$ or $M_2$, respectively.
		\begin{itemize}
			\item[\textnormal{(i)}] A $\Gamma$-\emph{module} \emph{homomorphism} is a monoid homomorphism $\phi:M_1\rightarrow M_2$ that respects the action of $\Gamma$, this means $\phi(^\alpha a)={^\alpha \phi(a)}$.
			\item[\textnormal{(ii)}] A $\Gamma$-\emph{order-ideal} of a monoid $M$ is a subset $I$ of $M$ such that for any $\alpha, \beta\in \Gamma$, $^\alpha a+{^\beta b}\in I$ if and only if $a,b\in I$.
\item[\textnormal{(iii)}] $M$ is said to be a $simple$ $\Gamma$-monoid if the only $\Gamma$-order-ideals of $M$ are $0$ and $M$.
	\end{itemize}

\end{defn}

\begin{prop}\textnormal{\cite{A2}} 
\label{1prop3}
Let $\Gamma$ be a group and $T$ a $\Gamma$-monoid. For $x\in T$, the $\Gamma$-order-ideal  generated by $x$, denoted by $\langle x\rangle$, is given by the set
$$ \langle x\rangle =\left \{ y\in T:y\leq \sum_{\alpha \in \Gamma} k_\alpha{^\alpha}x  \right \}.$$
\end{prop}


\begin{defn}
\textnormal{\cite{Alahmedi, A1,hazratGK, namzerui}} A (finite) path $p=\alpha_1\alpha_2\cdots \alpha_n$ such that  $s(\alpha_1)=r(\alpha_n)$ and $s(e_i)\neq s(e_j)$ for all $i\neq j$  is called a $cycle$. A cycle of length $1$ is called a \emph{loop}. Two cycles $C$ and $D$ is called \emph{distinct} if $C^0\neq D^0$. An edge $f$ is an \emph{exit} for a path $p=e_1\cdots e_n$ if $s(f)=s(e_i)$ but $f\neq e_i$ for some $1\leq i\leq n$. A graph $E$ is said to be a \emph{no-exit graph} if no cycle in $E$ has an exit. Two cycles are said to be \emph{disjoint} if they do not have any common vertices. The set of all paths in $E$ shall be denoted by Path$(E)$. For nonempty subsets $X$ and $Y$ of $E^0$, we denote $E(X,Y)=\{ e\in E^1: s(e)\in X, r(e)\in Y \}$. 
\end{defn}



\begin{defn} \textnormal{\cite{A1_8}}
For any vertex $v$, the \emph{tree} of $v$, denoted by $T(v)$, is the set $\{w\in E^0: v=s(p) \textnormal{~and~} w=r(p), p\in Path(E)  \}$. A finite graph $E$ is called a $C_n$-comet, if $E$ has exactly one cycle $C$ (of length $n$), and $T(v)\cap C^0\neq \varnothing $ for any vertex $v\in E^0$. $E$ is called a
\emph{multi-headed comet} if $E$ consists of $C_{l_s}$-comets, $1\leq s\leq t$, of length $l_s$, such that the cycles are mutually disjoint and for any $v\in E^0$, there exists $ l_k\leq t$ such that $T(v)\cap (C_{l_k})^0\neq \varnothing$, and furthermore, $E$ is a no-exit graph.
\end{defn}

\begin{defn} \cite{LPAbook}
A graph $E$ is \emph{connected} if $\widehat{E}$
is a connected graph in the usual sense, that is, if given any two vertices $u,v\in E^0$, there exist $h_1,h_2,\cdots ,h_m\in E^1\cup (E^1)^*$ such that $p=h_1h_2\cdots h_m$ is a path in $\widehat{E}$ with $s(p)=u$ and $r(p)=v$. This is equivalent to saying that the underlying undirected graph of $E$ is connected. 
\end{defn}


\begin{defn} {\textnormal{\cite{A1}} Let $E$ be a row-finite directed graph. The \emph{graph monoid} of $E$, denoted by $M_E$, is the abelian monoid generated by $\lbrace v:v\in E^0\rbrace $, subject to
$$v= \sum_{e\in s^{-1}(v) }r(e)$$ for every $v\in E^0$ that is not a sink.

 The \emph{talented monoid} of $E$, denoted by $T_E$, is the abelian monoid generated by $\lbrace v(i): v\in E^0, i \in \mathbb{Z}\rbrace$, subject to 
$$v(i)= \sum_{e\in s^{-1}(v) }r(e)(i+1)$$
for every $i \in \mathbb{Z}$ and every $v\in E^0$ that is not a sink. 
 }\end{defn}

The talented monoid $T_E$ is equipped by a natural $\mathbb{Z}$-action:
$$^nv(i)=v(i+n), $$
which implies the talented monoid $T_E$ is a $\mathbb{Z}$-monoid. 

\noindent Throughout the article, we simultaneously use $v\in E^0$ as a vertex in $E$, as an element of $L_K(E)$, and the element $v= v(0)$ in $T_E$. Hence, $^n v ={^n}v(0)= v(n)$. In the same manner, for  $w\in T_E$, by saying $w\in E^0$ simply means $w=w(0)$.

\begin{defn}{\cite{A1}
 Let $\Gamma$ be an abelian group with identity $e$. A ring $A$ (possibly without unit)
is called a $\Gamma$-\emph{graded ring} if 
$A=\oplus_{\gamma\in \Gamma}A_\gamma$ such that each $A_\gamma$ is an additive subgroup of $A$ and 
$A_\gamma A_\delta\subseteq A_{\gamma \delta}$ for all $\gamma,\delta \in \Gamma$. The group $A_\gamma$ is called a $\gamma$-$homogenous$ $component$ of $A$.  The elements of $\bigcup_{\gamma \in \Gamma} A_\gamma$ are called $homogenous $ $ elements$ of $A$. The nonzero elements of $A_\gamma$ are called $homogenous$ $of$ $degree$ $\gamma$ and we write deg$(a)=\gamma$ for $a\in A_\gamma\setminus \{0\}$. When it is clear from the context that a ring $A$ is graded by the group $\Gamma$, we simply say that $A$ is a graded ring.

 An ideal $I$ of $A$ is called a \emph{graded ideal} if it is generated by homogenous elements. This is equivalent to
$I=\oplus_{\gamma \in \Gamma }I\cap A_\gamma.$ We say that $A$ is a $graded$ $simple$ $ring$ if the only ideals of $A$ are $0$ and $A$.}
\end{defn}


 Setting deg$(v)=0$ for $v\in E^0$, deg$(\alpha)=1$ and deg$(\alpha^*)=-1$ for $\alpha \in E^1$, we obtain a natural $\mathbb{Z}$-grading on the free $K$-ring generated by $\{v,\alpha , \alpha^* :v\in E^0, \alpha \in E^1  \}$. Since the relations in definition the Leavitt path algebra are all homogenous, the ideal generated by these relations is homogenous and thus we have a natural $\mathbb{Z}$-grading on $L_K(E)$.

Let $L_K(E)$ be a Leavitt path algebra with coefficients in the field $K$ associated to the row-finite  graph $E$. $E$ is said to be \emph{simple} if the Leavitt path algebra $L_K(E)$ is  simple. Denote by $\mathcal{L}^{gr}(L_K(E))$ the lattice of graded ideals of $L_K(E)$, $\mathcal{L}(E)$ the set of hereditary saturated subsets of $E$, and $\mathcal{L}(T_E)$ the lattice of $\mathbb{Z}$-order-ideals of $T_E$. For $H\subseteq E^0$, let $I(H)$ and $\langle H\rangle $  be the graded ideal and the $\mathbb{Z}$-order ideal generated by $H$, respectively. For $v\in E^0$, the $\mathbb{Z}$-order-ideal of $T_E$ generated by $\{v\}$ shall be denoted by $\langle v\rangle$. \\
\indent The following results provides an overview of the relationship of the geometry, the monoid structure and the algebraic structure associated to a graph.

\begin{thm}{\textnormal{ \cite{LPAbook}}}\label{thm29} 
Let $E$ be a row-finite graph. Then there is a lattice isomorphism between $\mathcal{L}(E)$ and $\mathcal{L}^{gr}(L_K(E))$ [resp., between $\mathcal{L}(E)$ and $\mathcal{L}(T_E)$, between $\mathcal{L}(T_E)$ and $\mathcal{L}^{gr}(L_F(E))$] given by 
$H\mapsto I(H)$ [resp., $H\mapsto \langle H\rangle $, $\langle H\rangle \mapsto I(H)$]
where $H$ is a hereditary saturated subset of $E$.

\end{thm}

It then follows that every $\mathbb{Z}$-order-ideal of the talented monoid is generated by some hereditary saturated set, which essentially composes the set of vertices in the ideal.

Since then, algebraic results on the Leavitt path algebras were ``translated" into its monoid counterparts. One of the many attempts made to solidify this relationship is the investigation of the  Gelfand-Kirillov dimension of the Leavitt path algebras which was translated in the language of talented monoids in \cite{hazratGK}.

\begin{defn}{ \textnormal{\cite{GKdimension}} Let $A$ be an algebra (not necessarily unital), which is generated by a finite dimensional subspace $V$. Let $V^n$ denote the span of all products $v_1v_2\cdots v_n,$
$v_i\in V$, $k\leq n$. Then $V=V^1\subseteq V^2\subseteq \cdots$,
\begin{center}
$A=\bigcup_{n\geq 1}V^n~~$
and $~~g_{V(n)}=\textnormal{dim}V^n<\infty$.
\end{center}
Given the functions $f,g:\mathbb{N} \rightarrow \mathbb{R}^+$, if there exists $c\in \mathbb{N} $ such that $f(n)\leq cg(cn)$ for all $n\in \mathbb{N} $ we call $f$ asymptotically bounded by $g$.\\
\indent If $f$ is asymptotically bounded by $g$ and $g$ is asymptotically bounded by $f$, the functions $f$ and $g$ are said to be $asymptotically$ $ equivalent$ denoted by $f\sim g$. The equivalence class of $f$ under $\sim$ is called the $growth$ of $f$.\\
\indent If $W$ is another finite-dimensional subspace that generated $A$, then $g_{V(n)}\sim g_{W(n)}$. If $g_{(V(n))}$ is polynomially bounded, then  the \emph{Gelfand-Kirillov dimension} (GK-dimension) of $A$ is defined as 
\begin{center}
GKdim$A$ $=$ $\displaystyle \limsup_{n\rightarrow \infty } \dfrac{\ln g_{V(n)}}{\ln n}.$
\end{center}
The GK-dimension does not depend on a choice of the generating space $V$ as long as dim$V<\infty$. If the growth of $A$ is not polynomially bounded, then GKdim$A=\infty$.
}
\end{defn}

\begin{thm}\textnormal{\cite{zel}}\label{zelgk} 
Let $E$ be a finite graph.

\begin{enumerate}
\item [\textnormal{(1)}] The Leavitt path algebra $L_K(E)$ has polynomially bounded growth if and only if $E$ is a graph with disjoint cycles. 

\smallskip

\item [\textnormal{(2)}] If $d_1$ is the maximal length of a chain of cycles in $E$, and $d_2$ is the maximal length of chain of cycles with an exit, then $$GKdim L_K(E) = \max(2d_1-1, 2d_2).$$
\end{enumerate}
\end{thm}

Throughout, we shall constantly call the formula for the GK-dimension found in Theorem \ref{zelgk}, and $d_1$ and $d_2$ shall  denote the maximal length of a chain of cycles in $E$, and maximal length of chain of cycles in $E$ which has an exit, respectively.

\begin{defn}\cite{sebandalvilela} \label{dfh1} {Let $I$ be a $\Gamma$-order-ideal of a $\Gamma$-monoid $T$. We say
	\begin{enumerate}
		\item[\textnormal{(i)}] $I$ is a $cyclic$ $ideal$ if for any $x\in I$, there is an $\alpha\in 
		\Gamma$ such that $^\alpha x=x$;
		\item[\textnormal{(ii)}] $I$ is a $comparable$ $ideal$ if for any $x\in I$, there is an $\alpha \in \Gamma$ such that $^\alpha x> x$;
		\item[\textnormal{(iii)}] $I$ is a $non$-$comparable$ $ideal $ if for any $x\in I$, and any $\alpha \in \Gamma$, we have $^\alpha x\parallel x$.
	\end{enumerate}
	}
\end{defn}

\begin{defn} \cite{sebandalvilela} \label{df45} {Let $T$ be a $\Gamma$-order-ideal. A \emph{$\Gamma$-series} for $T$ is a sequence of $\Gamma$-order-ideals $$~~~~~~~~~~~~~~~~~~~0=I_0 \subseteq I_1 \subseteq I_2 \subseteq \cdots \subseteq I_n=T.~~~~~~~~~~~~~~~~~~~(*) $$  The \emph{length} of a $\Gamma$-series is the number of its proper inclusions. A $refinement$ of $(*)$ is any $\Gamma$-series of the form 
		$$~~~~~~~~~~~~~~~~~~~0=I_0 \subseteq I_1 \subseteq \cdots \subseteq I_i\subseteq N\subseteq I_{i+1}\subseteq \cdots \subseteq  I_n=T,~~~~~~~~~~~~~~~~~~~ $$ and this refinement is said to be $proper$ if $I_i\subsetneq N\subsetneq I_{i+1}$.
		Furthermore, we say $(*)$ is a  $\Gamma$-$composition$ $series$ if for each  $i=0,1\cdots,n-1$, $I_i\subsetneq I_{i+1}$ and each of  quotients $I_{i+1}/I_{i}$ are simple $\Gamma$-monoids.\\ \indent  A $\Gamma$-composition series is of $cyclic$ [\textit{resp., non-comparable}, \textit{comparable}] type if all of the simple quotients $I_{i+1}/I_{i}$ are cyclic [resp., non-comparable, comparable]. 
	}
\end{defn}

For a graph $E$ with disjoint cycles, the composition series of the $\mathbb{Z}$-monoid $T_E$ was then characterized in relation to the Gelfand-Kirillov dimension.

\begin{lemma} \textnormal{\cite{hazratGK}} \label{hazratGK-4.1} Let $E$ be a finite graph and $T_E$ its talented monoid. Then the following are equivalent.
\begin{enumerate}
    \item [\textnormal{(i)}] $E$ consists of disjoint cycles with no sinks.
    \item [\textnormal{(ii)}] $T_E$ has a cyclic composition series. 
\end{enumerate}
\end{lemma}

\begin{thm}\label{mainthm} \textnormal{\cite{hazratGK}}
Let $E$ be a finite graph, $L_K(E)$ its associated Leavitt path algebra and $T_E$ its talented monoid.  Then the following are equivalent.

\begin{enumerate}
\item [\textnormal{(i)}] $E$ is a graph with disjoint 
 cycles.
\item [\textnormal{(ii)}] $T_E$ has a composition series of cyclic and non-comparable types. 
\item [\textnormal{(iii)}] $L_K(E)$ has a finite GK-dimension. 
\end{enumerate}

\end{thm}

\section{No-exit Graphs}

We present some properties  no-exit graphs relating to the ideal structure of the talented monoid and the GK-dimension of the associated Leavitt path algebra. 

\begin{lemma}\label{2lem1} 
$E$ is a finite no-exit graph if and only if  $GKdim~L_K(E)\leq 1$.
\end{lemma}

\pf Suppose $E$ is a no-exit graph, that is, every cycle in $E$ has no exit. Consider the following cases.

\indent Case 1. $E$ has no cycles. 

\indent Then length of chain of cycles in $E$ must be $0$. Thus, $d_1=0=d_2$. By Theorem \ref{zelgk}, $GKdim L_K(E) = \max(2d_1-1, 2d_2) = \max(2(0)-1, 2(0))=\max(-1,0)=0$.

\indent Case 2. $E$ has a cycle. 

\indent Since every cycle in $E$ has no exit, and $E$ has at least one cycle, the maximum length of a chain of cycles in $E$ $d_1$ must be $1$. It also follows that no chain of cycles in $E$ has an exit. Thus, $d_2=0$. Accordingly, $GKdim L_K(E) = \max(2d_1-1, 2d_2)=\max(2(1)-1, 2(0))=1$.

Conversely, suppose $GKdim L_K(E)=\max(2d_1-1, 2d_2)\leq 1$.  Then $2d_1-1\leq 1$ and $2d_2\leq 1$. That is, $d_1\leq 1$ and $d_2\leq \frac{1}{2}$. Since the GK-dimension is a nonnegative integer, $d_2$ must be $0$. Since $d_2$ is the maximal length of a chain of cycles with an exit and $d_2=0$, it follows that no cycle in $E$ has an exit. Hence, $E$ is a no-exit graph. \qed

\begin{prop}
\label{2lem5}
Let $E$ be a row-finite graph. Then the following are equivalent
\begin{enumerate}
\item [\textnormal{(i)}] for every $v\in E^0$, $\langle v\rangle \cap E^0=\{v\}$;
\item [\textnormal{(ii)}] for any distinct $v,w\in E^0$, $v\not \in r(s^{-1}(w))$.
\end{enumerate}

\end{prop}

\pf 
Suppose for every $v\in E^0$, $\langle v\rangle \cap E^0=\{v\}$. In contrary, suppose there exists $v_0\neq w_0$ with $v_0\in r(s^{-1}(w_0))$. Then
$$^{-1}w_0 =\sum_{e\in s^{-1}(w_0)}r(e)(0)= \sum_{r(e)\neq v_0}r(e)(0)+v_0(0).$$
Thus, $^0v_0\leq {^{-1}}w_0$ which implies $v_0\in \langle w_0\rangle$ with $v_0\neq w_0$. Accordingly, $\langle w_0\rangle \cap E^0\neq \{w_0\}$, a contradiction. Hence, for any distinct vertices $w$ and $v$, $v\not \in r(s^{-1}(w))$.\\

Conversely, suppose there exits $w\in E^0$ such that $\langle w\rangle\cap E^0\neq \{w\}.$ Then there exists $u\in E^0$, $u\neq w$ such that $u\in\langle w\rangle$. Thus, by Proposition \ref{1prop3},
$$u\leq \sum_{\alpha \in \Gamma}k_\alpha {^\alpha} w.$$
Hence, there is a path from $w$ to $u$, that is, there exists $p=e_1\cdots e_n\in \textnormal{Path}(E)$ such that $s(p)=w$ and $r(p)=u$. Since $u\neq w$, there exists $m\leq n$ such that $r(e_m)\neq w$. Suppose $m'$ is the least such $m$ with $r(e_m)\neq w$. Hence, $p'=e_1 e_2\cdots  e_{m}\in \textnormal{Path}(E)$ with $s(p')=w$ and $r(p')=r(e_m)\neq w$. Accordingly, $r(e_m) \in r(s^{-1}(w))$ with $r(e_{m})\neq w$, a contradiction. \qed

\begin{rema} \label{6rem1}
Notice that statement $(ii)$ of Proposition \ref{2lem5} is equivalent to saying that every vertex in $E$ are not connected to each other by an edge. Hence, it follows that  $\langle v\rangle\cap E^0=\{v\}$ for every $v\in E^0$ if and only if for any $w\neq v$, $v$ and $w$ are not connected by an edge. 
\end{rema}

\begin{rema}
If a finite graph $E$ is a disjoint union of isolated vertices and loops, then $E$ is a graph with disjoint cycles. By Theorem \ref{mainthm}, it follows that $T_E$ has a composition series of cyclic and non-comparable types and $L_K(E)$ has a finite GK-dimension. 
\end{rema}

\noindent The following proposition particularly computes the GK-dimension of graphs which is a disjoint union of isolated vertices and loops and  provides a characterization for such graphs.

\begin{prop}
\label{2lem6}  Let $E$ be a finite graph. Then $E$ is a disjoint union of isolated vertices and loops if and only if $GKdim~L_K(E)\leq 1$ and for every $v\in E^0$, $\langle v\rangle \cap E^0=\{v\}$.
\end{prop}

\pf Suppose $E$ is a disjoint union of vertices and loops. Then clearly, $E$ is a no-exit graph. By Lemma \ref{2lem1}, $GKdim~L_K(E)\leq 1$. \\
\indent Suppose there exists $w\in E^0$ such that $\langle w\rangle\cap E^0\neq \{w\}.$ Then there exists $u\in E^0$, $u\neq w$ such that $u\in\langle w\rangle$. Thus, by Proposition \ref{1prop3},
$$u\leq \sum_{\alpha \in \Gamma}k_\alpha {^\alpha} w.$$
Hence, there is a path from $w$ to $u$. Since $u\neq w$, $E$ cannot be a disjoint union of isolated vertices and cycles, a contradiction. 
Hence, $\langle v\rangle \cap E^0=\{v\}$ for every $v\in E^0$.\\
\indent Conversely, suppose $GKdim~L_K(E)\leq 1$ and for any $v\in E^0$, we have $\langle v\rangle \cap E^0=\{v\}$. 
\indent Assume $GKdim~L_K(E)= 0$. Then $\max(2d_1-1, 2d_2)=0$, that is, $2d_1-1\leq 0$ and $2d_2\leq 0$. We then obtain $d_1\leq \frac{1}{2}$ and $d_2\leq 0$. Since $d_1$ and $d_2$ are nonnegative integers, we have $d_1=0=d_2$. Hence, $E$ must not have any cycles. In particular, $E$ has no loops. \\
\indent By Remark \ref{6rem1}, no two distinct vertices are connected by an edge. Since this case is only true for isolated vertices and loops and $E$ has no loops, it follows that $E$ must be composed of isolated vertices. \\
\indent Assume $GKdim~L_K(E)=1$. Then $\max(2d_1-1,2d_2)=1$, that is, $2d_1-1\leq 1$ and $2d_2\leq 1$. We then have $d_1\leq 1$ and $d_2\leq \frac{1}{2}$, which also implies $d_2=0$. \\
\indent If $d_1=0$, then by previous arguments, it follows that $E$ is a disjoint union of isolated vertices. \\
\indent Suppose $d_1=1$. Then $E$ has at least one cycle and since $d_2=0$, this cycle has no exit. Now, let $C$ be a cycle in $E$. By Proposition \ref{2lem5}, for any  distinct vertices $v,w\in E^0$, $v\not \in r(s^{-1}(w))$, that is, no two distinct vertices are connected by an edge. This implies that $C$ must be a loop. Since $GKdim~L_K(E)=1<\infty$, by Theorem \ref{mainthm}, cycles in $E$ must be disjoint, that is in particular, $C$ is disjoint to any other loops. \\
\indent Therefore, $E$ is a disjoint union of isolated vertices and loops. \qed

\begin{thm}
\label{2lem4} 
Let $E$ be a row-finite no-exit graph and $v\in E^0$. Then the following are equivalent for a vertex $v$:
\begin{enumerate}
    \item [\textnormal{(i)}]
    for every $e\in s^{-1}(v)$, $r(e) \parallel w$ for all $w\neq v$ as elements of $M_E$ (in particular, $v=\sum_{e\in s^{-1}(v)}{^1}r(e)$ is a unique representation of $v$ (by indexing)); 
\item [\textnormal{(ii)}] for each $e\in s^{-1}(v)$, $r(e)$ is either a sink with $r^{-1}(r(e))=\{e\}$, or in a loop $f$ with $r^{-1}(r(e))=\{e,f\}$. 
\end{enumerate}
\end{thm}

\pf $(\Rightarrow)$ Suppose there exists $e_0\in s^{-1}(v)$ such that $r(e_0)$ is not a sink  or $r^{-1}(r(e_0))\neq \{e_0\}$. 
Then $s^{-1}(r(e_0))\neq \varnothing$  or $\{e_0\} \subsetneq r^{-1}(r(e_0))$,  since $e_0\in r^{-1}(r(e_0))$. 
\\
\indent Suppose $s^{-1}(r(e_0))\neq \varnothing$. Let $d\in s^{-1}(r(e_0))$. If $r(e_0)\neq r(d)$, then, by our assumption, $r(e_0)\parallel r(d)$. Thus, $r(e_0)\not >r(d)$. This in particular implies that no edge must connect $r(e_0)$ to $r(d)$, a contradiction since $d\in s^{-1}(r(e_0))$. Hence, $r(e_0)=r(d)$, that is, $d\in r^{-1}(r(e_0))$ and $d$ is a loop. \\
\indent If there exists $u\neq v$ such that $r(e_0)=r(e_u)$ for some $e_u\in s^{-1}(u)$, then $u>r(e_0)$,  contradiction. Hence, $r^{-1}(r(e_0))=\{e_0,d\}$. \\
\indent Suppose $\{e_0\}\subsetneq r^{-1}(r(e_0))$. Then there exists $g\neq e_0$ such that $g\in r^{-1}(r(e_0))$. Then $r(g)=r(e_0)$. Suppose $s(g)\neq r(e_0)$. Then $s(g)$ is connected to $r(g)=r(e_0)$ implying that $s(g)>r(e_0)$. That is, $s(g)\nparallel r(e_0)$, a contradiction. Hence, $s(g)=r(e_0)=r(g)$, that is, $g$ is a loop. Similarly above, we must also have  $r^{-1}(r(e_0))=\{e_0,g\}$. \\

$(\Leftarrow)$ Let $v\in E^0$. Then
$$v=\sum_{e\in s^{-1}(v)}r(e)(1).$$
\indent Let $A\subseteq s^{-1}(v)$ be the set of all edges $e$ in $s^{-1}(v)$ such that $r(e)$ is a sink with $r^{-1}(r(e))=\{e\}$ and let $B$ be the set of all edges $e$ in $s^{-1}(v)$ such that $r(e)$ is in a loop $f_e$ with $r^{-1}(r(e))=\{e,f_e\}$. Then
$$v(i)=\sum_{e\in A}r(e)(i+1) + \sum_{e\in B}r(e)(i+1) .$$

Let $e\in A$. Then as a sink, $s^{-1}(r(e))=\varnothing$ and since $r^{-1}(r(e))=\{e\}$, there exists no vertex $w\neq v$ such that $w$ and $v$ are connected by an edge. Hence, $w\parallel v$ for all $w\neq v$ as elements of $M_E$ and $r(e)$ cannot be further transformed.\\

Let $e\in B$. Then $e$ is in a loop $f_e$ with $r^{-1}(r(e))=\{e,f_e\}$. Since $E$ is a no-exit graph, $f_e$ has no exits, that is, $s^{-1}(r(e))=\{f\}$. Thus, 
$$ \hspace{2cm} r(e)=r(e)(1)=r(e)(2)=\cdots~.  \hspace{2cm} (*)$$
Also, since $r^{-1}(r(e))=\{e,f\}$, $s(e)=v$ is the only vertex connected to $r(e)$. Thus, $(*)$ is the only transformation for $r(e)$.
Accordingly, we obtain 
$$v(i)=\sum_{e\in A}r(e)(i+1) + \sum_{e\in B}r(e)(1) .$$
That is, for $i=0$,
$$\hspace{4.1cm}v=v(0)=\sum_{e\in s^{-1}(v)}r(e)(1).\hspace{4.1cm}\qed$$


\section{Solvability and Nilpotency}

\begin{defn}\textnormal{\cite{namzerui}
Let $(L, [-,-])$ be a Lie algebra. Define 
\begin{center}$L^{(0)}= L = L^0$,~ 
$L^{(n)}=[L^{(n-1)}, L^{(n-1)}]$~
and ~$L^{n}=[L, L^{n-1}]$
\end{center} for every $n\geq 1$.
Then $L$ is called \emph{solvable} (resp. \emph{nilpotent}) of index $n$ if $n$ is the minimal integer such that $L^{(n)}=0$ (resp. $L^n=0$)}.
\end{defn}

\begin{defn}\textnormal{\cite{namzerui}
Let $K$ be a field and $A$ an associative $K$-algebra. Then $A$ becomes a Lie algebra under the operation $[x,y]=xy-yx$ for all $x,y\in A$, and $(A, [-,-])$ is called the \emph{associated Lie algebra} of $A$. \\
\indent The associative algebra $A$ is called \emph{Lie solvable} (resp. \emph{Lie nilpotent}) of index $n$ if its associated Lie algebra is solvable (resp. nilpotent) of index $n$. }
\end{defn}

\begin{lemma} \label{1lem1}\textnormal{\cite{namzerui}} 
Let $K$ be an arbitrary field and 
 $E$ an arbitrary graph such that $L_K(E)$ is Lie solvable. Then $E$ is a no-exit graph.
\end{lemma}

\begin{lemma}\label{2lem2} 
Let $K$ be an arbitrary field and 
 $E$ a finite graph such that $L_K(E)$ is Lie solvable. Then $GKdim~L_K(E)\leq 1$. 
\end{lemma}
\pf Suppose $L_K(E)$ is Lie solvable. By Lemma \ref{1lem1}, $E$ is a no exit graph. By Lemma \ref{2lem1}, $GKdim~L_K(E)\leq 1$. 
\qed
~\\~\\
Notice that in a row-finite no-exit graph $E$, no two cycles would have a common vertex. Otherwise, if $C_1$ and $C_2$ have a vertex $v$ in common, then $e\in C_2^1$ with $s(e)=v$ would be an exit for the cycle $C_1$. Hence, $E$ is a graph with disjoint cycles. By Theorem \ref{mainthm} and Lemma  \ref{2lem2}, the connection by the following corollary could be seen. 

\begin{cor}
Let $K$ be an arbitrary field and $E$ a finite graph such that $L_K(E)$ is Lie solvable. Then $E$ is a graph with disjoint cycles and has a composition series of cyclic and noncomparable types. 
\end{cor}

\begin{thm} \label{1thm2} \textnormal{\cite{namzerui}} 
Let $K$ be an arbitrary field and $E$ an arbitrary graph. Then the following hold:
\begin{enumerate}
    \item [\textnormal{(i)}] If $\textnormal{char}(K)=2$, then $L_K(E)$ is Lie solvable if and only if $E$ is a no-exit graph satisfying the following condition: every vertex $v\in E$ is a sink, or in a cycle whose length is at most 2, or for each $e\in s^{-1}(v)$, $r(e)$ is either a sink with $r^{-1}(r(e))=\{e\}$, or in a loop $f$ with $r^{-1}(r(e))=\{e,f\}$. 

\item [\textnormal{(ii)}] If $\textnormal{char}(K)\neq 2$, then $L_K(E)$ is Lie solvable if and only if $E$ is a disjoint union of isolated vertices and loops. In this case, $L_K(E)\cong K^{(\gamma_1)} \oplus K[x, x^{-1}]^{(\gamma_2)}$, where $\gamma_1$ is the set of all isolated vertices in $E$ and $\gamma_2$ is the set of all loops in $E$. 
\end{enumerate}
\end{thm}

\begin{cor}\label{2rem3}
Let $K$ be a field with char$(K)=2$ and $E$ be a row-finite graph with $L_K(E)$ Lie solvable. Then  $|\langle v\rangle \cap E^0|=2$ for any $v$ in a cycle $C$ if and only if $C$ is of length $2$. 
\end{cor}
 
\pf 
 Let $K$ be a field with char$(K)=2$ and $E$ be a row-finite graph with $L_K(E)$ Lie solvable. If $E$ is has a cycle $C$, then by Theorem \ref{1thm2}(i), it must be  of length at most $2$ and for every $v\in E^0\setminus C^0$, we have $r(s^{-1}(v))\cap C^0=\varnothing$. That is, no edge has range in $C$. In contrary, suppose there exists  $u\in r(s^{-1}(v))\cap C^0$. Then for $e\in r^{-1}(u)\cap s^{-1}(v)$, $s(e)=v$ is not a sink. Since $E$ is a no-exit graph, $s(e)=v$ is not in a cycle of length $2$. Also,  $r(e)=u$ is neither a sink nor in a loop (otherwise, a loop in $u$ is an exit in $C$). These all imply that a cycle in $E$ must be isolated. Hence, the number of vertices in $\langle v\rangle $ is the length of the cycle which is $2$. The converse is clear to see. \qed

\begin{thm}\textnormal{\cite{hazratGK}} \label{hazratGK-3.11} Let $E$ be a finite graph. There is a one-to-one correspondence between sinks in $E$ and the non-comparable minimal ideals of $T_E$. 
\end{thm}

 The following theorem is a direct consequence of Theorems \ref{1thm2} and \ref{2lem4}, and Lemmas  \ref{2lem1}, \ref{2rem3}, \ref{2lem6} and \ref{hazratGK-3.11}.

\begin{thm}
\label{2thm7}  Let $K$ be an arbitrary field and $E$ a finite graph. Then the following hold:
\begin{enumerate}
    \item [\textnormal{(i)}] In the case $\textnormal{char}(K)=2$, $L_K(E)$ is Lie solvable if and only if we have $GKdim~L_K(E)\leq 1$ and one of the following conditions is satisfied: for every vertex $v\in E$, we have either $\langle v\rangle$ is a minimal non-comparable ideal, or a cyclic ideal with $|\langle v\rangle \cap E^0|\leq 2$, or for every $e\in s^{-1}(v),$ $r(e) \parallel w$ as elements of $M_E$  for every $w\neq v$.

\item [\textnormal{(ii)}] For $\textnormal{char}(K)\neq 2$, we have $L_K(E)$ is Lie solvable if and only if we have $GKdim~L_K(E)\leq 1$ and for every $v\in E^0$, $\langle v\rangle \cap E^0=\{v\}$. 
\end{enumerate}
\end{thm}

\begin{cor} \label{1cor6} \textnormal{\cite{namzerui}}
For every field $K$ and every graph $E$, the following conditions are equivalent:
\begin{enumerate}
    \item [\textnormal{(i)}] $L_K(E)$ is Lie nilpotent;
    \item [\textnormal{(ii)}]  $E$ is a disjoint union of isolated vertices and loops;
    \item [\textnormal{(iii)}] $L_K(E)\cong K^{(\gamma_1)}\oplus K[x,x^{-1}]^{(\gamma_2)}$, where $\gamma_1$ is the set of all isolated vertices in $E$ and $\gamma_2$ is the set of all loops in $E$. 
\end{enumerate}
\end{cor}

\begin{rema}
By Theorem \ref{1thm2} and Corollary \ref{1cor6}, it then follows that for fields of characteristic not equal to $2$, Lie solvability and Lie nilpotency are equivalent. 

\end{rema}
The following result gives a characterization of Lie nilpotency in terms of the GK-dimension and the ideal structure of the talented monoid which directly follows from Corollary \ref{1cor6} and Lemma \ref{2lem6}.

\begin{thm} \label{2thm8} 
For every field $K$ and every finite graph $E$, $L_K(E)$ is Lie nilpotent if and only if $GKdim~L_K(E)\leq 1$ and for every $v\in E^0$, $\langle v\rangle \cap E^0=\{v\}$.
\end{thm}


\begin{cor} \label{1cor8}  \textnormal{\cite{namzerui}} 
For a field $K$ and a graph $E$, $L_K(E)$ is Lie solvable if and only if $[L_K(E), L_K(E)]$ is Lie nilpotent. 
\end{cor}
The following corollaries follow from Lemma  \ref{2lem2}, Corollary \ref{1cor8}, and Theorem \ref{2thm7}(ii).

\begin{cor}\label{2cor9} 
For a field $K$ and a finite graph $E$, if $[L_K(E), L_K(E)]$ is Lie nilpotent then $GKdim~L_K(E)\leq 1$.

\end{cor}

\begin{cor}\label{2cor10} 
Let $K$ be a field with $\textnormal{char}(K)\neq 2$ and $E$ a finite graph. We have the following equivalence:
\begin{itemize}
	\item[(i)] $[L_K(E),L_K(E)]$ is Lie nilpotent 
	\item[(ii)] $GKdim~L_K(E)\leq 1$ and  for every $v\in E^0$, $\langle v\rangle \cap E^0=\{v\}$. 
	\end{itemize}
\end{cor}



\section{Simplicity}
\noindent To understand the simplicity of the corresponding Lie bracket algebras, it is a first step of interest to take a look at the simplicity of Leavitt path algebras. Indeed in \cite{abramsmesyan} the Leavitt path algebra is always assumed to be simple. A characterization of simplicity for Leavitt path algebras can be found in e.g.~\cite{LPAbook,abramsmesyan}. 
\begin{thm} \label{simplicityLPA} \textnormal{\cite{abramsmesyan}} 
Let $K$ be a field, and let $E$ be a row-finite graph. Then $L_K (E)$ is simple if and only if $E$ has the following properties. 
\begin{enumerate}
    \item [\textnormal{(i)}] Every vertex $v$ of $E$ connects to every sink and every infinite path of $E$.
    \item [\textnormal{(ii)}] Every cycle of $E$ has an exit.
\end{enumerate}
In particular, if $E$ is finite, then $L_K (E)$ is simple if and only if every vertex $v$ of $E$ connects to every sink and every cycle of $E$, and every cycle of $E$ has an exit.
\end{thm}

The translation into the framework of talented monoids was already done in \cite{A1}.

\begin{prop} \textnormal{\cite{A1}} \label{A1-5.1}
Let $E$ be a row-finite graph and $L_K (E)$ its associated Leavitt path algebra.

\begin{enumerate}
    \item [\textnormal{(i)}] The algebra $L_K (E)$ is graded simple if and only if $T_E$ is simple.
     \item [\textnormal{(ii)}] The algebra $L_K (E)$ is simple if and only if $T_E$ is simple and for any $a \in T_E$, if $^n a$ and $a$ are comparable, $n \in \mathbb{Z}_{<0},$ then $^n a >a$.
\end{enumerate}
\end{prop}

The next corollary follows from Proposition \ref{A1-5.1}(ii) and \cite{abramsmesyan} Corollary 21.

\begin{cor}
Let $K$ be a field, and let $E$ be a row-finite graph having infinitely many vertices. If $T_E$ is simple and for any $a \in T_E$, if $^n a$ and $a$ are comparable, $n \in \mathbb{Z}_{<0},$ then $^n a >a$, we have that $[L_K (E), L_K(E)]$ is a simple Lie algebra.  
\end{cor}

\begin{lemma} \textnormal{\cite{Alahmedi}} \label{Alahmedi-lemma1} $[L_K (E) , L_K(E)] =0$ if and only if $E$ is a disjoint union of vertices and loops. 
\end{lemma}

The following remarks are direct consequences of Lemma \ref{Alahmedi-lemma1} together with Lemma \ref{2lem6}, and  Theorem \ref{1thm2} for the case characteristic $\neq 2$, respectively.

\begin{rema}
For a finite graph $E$, we have $[L_K (E) , L_K(E)] =0$ if and only if GKdim$L_K(E)\leq 1$ and for every $v\in E^0$, $\langle v\rangle \cap E^0=\{v\}$.
\end{rema}

\begin{rema}
For a finite graph $E$, and a field $K$ with  $char(K)\neq 2$, $L_K(E)$ is Lie solvable if and only if $[L_K (E) , L_K(E)] =0$. Hence, for $char(K)\neq 2$, $L_K(E)$ is Lie solvable only in the case of index $0$.
\end{rema}

In \cite{Alahmedi}, the concept of balloons were used used to characterize Lie simplicity of Leavitt path algebras. In Theorem \ref{4thm3}, we give a talented monoid perspective of balloons.

\begin{defn} \textnormal{\cite{Alahmedi}} \label{defballoon} We call a vertex $v$ in a connected graph $E$ a \emph{balloon} over a nonempty set $W\subseteq E^0$ if
\begin{enumerate}
    \item [\textnormal{(i)}] $v\not \in W$
    \item [\textnormal{(ii)}] there is a loop $C\in E(v,v)$
    \item [\textnormal{(iii)}] $E(v,W)\neq \varnothing$
    \item [\textnormal{(iv)}] $E(v, E^0)=\{C\} \cup E(v, W)$
    \item [\textnormal{(v)}] $E(E^0,v)=\{C\}$.
\end{enumerate}
\end{defn}

\begin{ex} Consider the graph $F$ below and let $W=\{t,u,z \}$. It is easy to see that $w$ is a balloon over $W$.   
\begin{figure}[H]
    \centering
    \includegraphics[width=5cm]{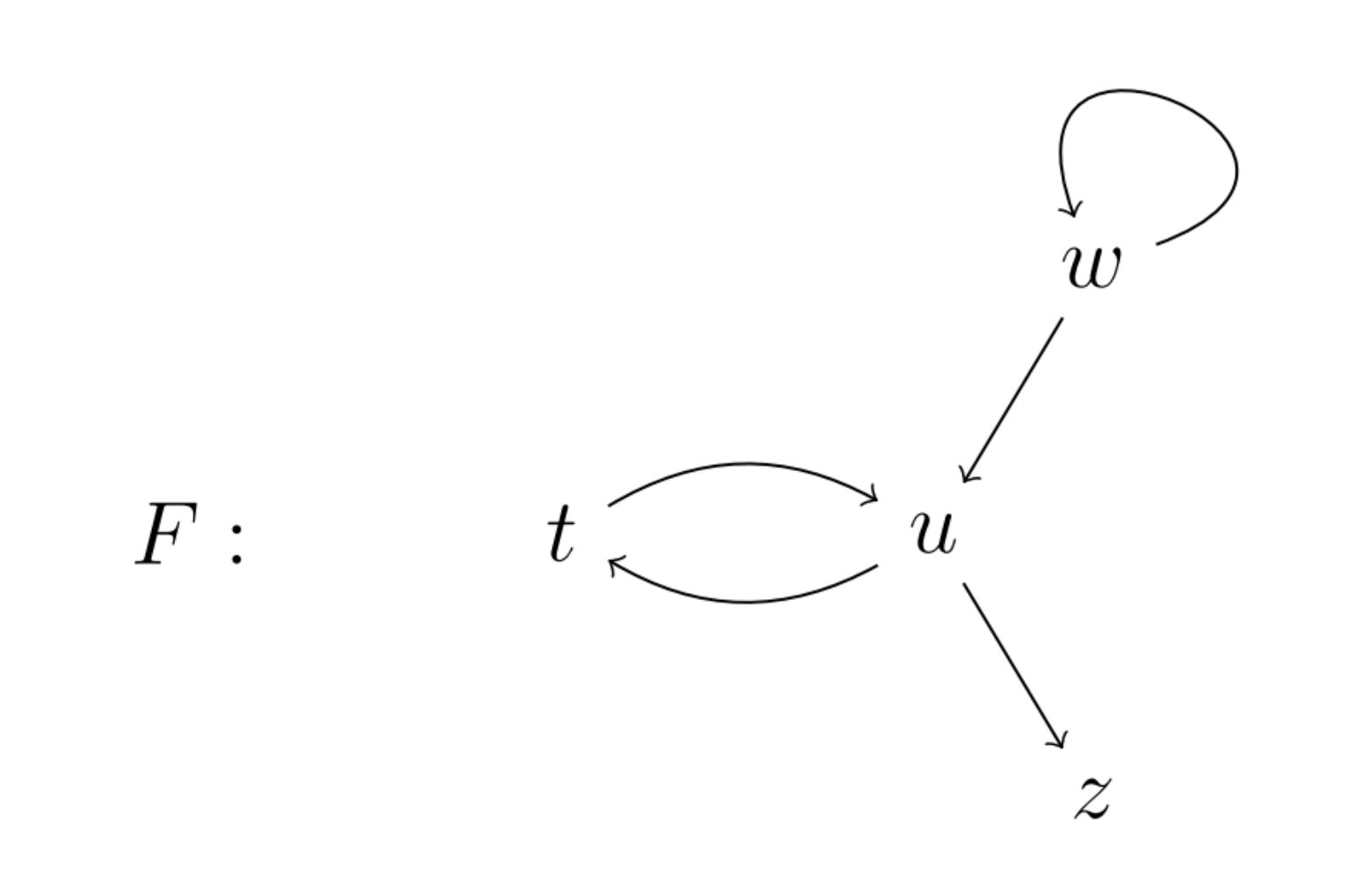}
    \label{balloon}
\end{figure}
\end{ex}

\begin{cor} \textnormal{\cite{hazratGK}} \label{hazratGK-3.12} 
Let E be a finite graph. Then the talented monoid $T_E$ is cyclic if and
only if $E$ is a multi-headed comet graph. Furthermore, $T_E$ is simple and cyclic if and
only if $E$ is a comet graph.

\end{cor}

\begin{thm} \label{4thm3}
 Let $E$ be a connected graph and $W\subseteq E^0$. A vertex $v\not \in  W$ is a balloon over $W$ if and only if 
\begin{enumerate}
    \item [\textnormal{(i)}] $H= E^0\setminus \{v\}$ is a hereditary saturated set (in that case, $\langle E\setminus \{v\}\rangle $ is the maximal $\mathbb{Z}$-order-ideal of $T_E$ which does not contain $v$);
    \item [\textnormal{(ii)}] 
    
    $r(s^{-1}(v))\setminus W =\{v\}$;

\item[\textnormal{(iii)}] $T_{E/H}$ is simple cyclic.
\end{enumerate}
\end{thm}

\pf Let $v$ be a balloon over $W$ and $H=E^0\setminus \{v\}$. Let $u\in H$ and $w\in r(s^{-1}(u))$. Since $E(E^0,v)=\{C\}$, where $C$ is a loop on $v$, $ v\not \in r(s^{-1}(u))$. Hence, in particular, $v\neq w $. Thus, $w\in H$, that is, $H$ is hereditary. Suppose $r(s^{-1}(u)) \subseteq H$ for some $u\in E^0$. Since $v$ has a loop $C$, $v\in r(s^{-1}(v))$, which implies that $r(s^{-1}(v))\not \subseteq H$. Thus, $u\neq v$, that is, $u\in H$. Hence, $H$ is saturated.


By the definition of a balloon, $E(v, E^0)=\{C\} \cup E(v,W)$. Let $e\in E(v,E^0)$. Since $x=r(e)\not \in W$, $e\notin E(v, W)$. Thus, $e=C$, that is $s(e)=v=r(e)=x$. Hence, $r(s^{-1}(v))\setminus W=\{v\}$. 

Notice that the quotient graph $E/H$ is an isolated loop, which is a comet graph. By Corollary \ref{hazratGK-3.12}, it follows that $T_{E/H}$ is simple cyclic.

Conversely, suppose $(i)-(iii)$ hold for some subset $W\subseteq E^0$ and a vertex $v\not \in W$. We show conditions $(ii)-(v)$ of Definition \ref{defballoon} hold. 



Now, by $(ii)$, $v\in r(s^{-1}(v))$. Thus, there is a loop $C\in E(v,v)$.


Assume that $E(E^0,v)\setminus\{C\}\neq \varnothing$. Then either there exists an edge  from $E^0\setminus \{v\}$ to $v$ or an edge $d$ aside from $C$ from $v$ to itself. In the first case, $H$ is not hereditary, which contradicts $(i)$. In the second case, it follows that $d$ is an exit for the cycle $C$. Hence, $C$ cannot be a comet graph. Thus, by Corollary \ref{hazratGK-3.12}, $T_{E/H}$ cannot be cyclic, a contradiction.

 Assume that there would be no edge from $v$ to $W$.
Since $E(E^0,v)=\{C\}$ as shown before and $E$ is connected we have that $E(v,W)\neq \varnothing$ and so, $E(v,E^0) = \{C\} \cup E(v,W)$.




 For a graph $E$, $H=E^0\setminus \{v\}$ being hereditary and saturated is equivalent to saying $\langle H\rangle$ is the largest $\mathbb{Z}$-order-ideal of $T_E$ which does not contain $v$. This completes the proof. 
\qed

    
    

\begin{rema}
Notice that $E/H$ is an isolated loop. By computation, we see that GKdim$L_K(E/H)=1$. By Theorem \ref{4thm3}(iii), $0\subseteq T_{E/H}$ is the only $\mathbb{Z}$-series for $T_{E/H}$. These coincide with Theorem \ref{mainthm}.
\end{rema}

\begin{lemma} \textnormal{\cite{A1}} \label{A1-2.2}
Let $E$ be a row-finite graph. For a hereditary saturated subset $H$ of $E^0$, and the $\mathbb{Z}$-order-ideal $I=\langle H\rangle\subseteq T_E$, we have the following $\mathbb{Z}$-module isomorphism
$$T_{E/H}\cong T_E/I.$$
\end{lemma}


\begin{thm} \textnormal{\cite{abramsmesyan}} \label{al1}  Let $E$ be a graph with $L_K(E)$ simple. 
\begin{enumerate}
    \item [\textnormal{(i)}] If $E^0$ is infinite, then $[L_K(E),L_K(E)]$ is simple.
    \item [\textnormal{(ii)}] If $E^0$ is finite, then $[L_K(E),L_K(E)]$ is simple if and only if 
    $$ 1_{L_K(E)}= \sum_{v\in E^0}v \not \in [L_K(E),L_K(E)]$$
\end{enumerate}
\end{thm}

\begin{thm}
\textnormal{\cite{Alahmedi}} \label{Alahmedi-Theorem2} Let $E$ be a row-finite graph. The Lie algebra $[L_K(E),L_K(E)]$ is simple if and only if either $L_K(E)$ is simple - in the case covered by Theorem \ref{al1} - or $E$ contains a simple subgraph $W$ such that every vertex $v\in E^0\setminus W$ is a balloon over $W$, and 
$$ \sum_{w\in r(E(v,W))}w \in [L_K(W), L_K(W) ]. $$
\end{thm}

The following results shall be a talented monoid viewpoint of the preceding results on Lie simplicity. 

\begin{thm}
Let $E$ be a row-finite graph where $L_K(E)$ is not simple and $[L_K(E), L_K(E)]$ is simple. Then $T_E/J$ has a cyclic composition series where $J$ is the minimal non-cyclic
$\mathbb{Z}$-order-ideal of $T_E$. 

\end{thm}

\pf Let $J=\langle W\rangle$, where $W$ is the intersection of all hereditary saturated sets in $E$. Then $J$ is a minimal $\mathbb{Z}$-order-ideal of $T_E$. By Theorem \ref{Alahmedi-Theorem2}, every vertex $v\not \in W$ is a balloon over $W$. By definition of a balloon, it follows that $E/W$ is a disjoint union of isolated loops. By Lemmas \ref{hazratGK-4.1} and \ref{A1-2.2}, $T_{E/W}=T_E/J$ has a cyclic composition series. 
\qed

It was shown in \cite{Alahmedi}, the set $W=\varnothing$ in Theorem \ref{Alahmedi-Theorem2} is the intersection of all hereditary saturated sets in $E$, which is of course, also a hereditary saturated set. $W$ then generates a minimal $\mathbb{Z}$-order-ideal of $T_E$. 
Thus, we have the following result.

\begin{thm}
\label{thmmain}
Let $E$ be connected row-finite graph with $L_K(E)$ not simple.  Then $[L_K(E), L_K(E)]$ is  simple if and only if  for every vertex $v\not \in I$ for some $\mathbb{Z}$-order-ideal $I$, (i)-(iii) of Theorem \ref{4thm3} are satisfied  
and 
$$ \sum_{w\in r(E(v,W))}w \in [L_K(W), L_K(W) ] $$
where $W=E^o\cap J$, $J$ the minimal non-cyclic $\mathbb{Z}$-order-ideal of $T_E$.

\end{thm}

\pf Let $v$ be a vertex not in some $\mathbb{Z}$-order-ideal $I$ of $T_E$. By Theorem \ref{thm29}, $I$ is generated by some hereditary saturated set in  $H= I\cap E^0$. Thus, $v\not \in H$, which  means 
$$v\not \in \bigcap \{ \textnormal{hereditary saturated sets in } E  \}=W,$$ 
by Theorem \ref{Alahmedi-Theorem2}. Hence, $v$ is a balloon over $W$, that is, $v$ satisfies (i)-(iii) of Theorem \ref{4thm3}. Now, being the intersection, it follows that $W$ generates a minimal $\mathbb{Z}$-order-ideal $J$ where $W$ precisely consists of the vertices in $J$. 

Suppose $J$ is cyclic. Then by Corollary \ref{hazratGK-3.12}, $W$ must be a comet graph. However, $W$ is a simple graph, which means $L_K(W)$ is simple. Thus, every cycle in $W$ must have an exit by Theorem \ref{simplicityLPA}, a contradiction to the definition of a comet graph.

Furthermore, by Theorem \ref{Alahmedi-Theorem2}, $$ \sum_{w\in r(E(v,W))}w \in [L_K(W), L_K(W) ]. $$This completes the proof. \qed

Now, Theorem \ref{4thm1} provides us conditions to where the Leavitt path algebra is simple. 

\begin{thm} \label{4thm1}
Let $E$ be a row-finite  graph. If $[L_K(E), L_K(E)]$ is simple and $L_K(E)$ is graded simple, then $L_K(E)$ is simple.
\end{thm}

\pf Suppose  $[L_K(E), L_K(E)]$ is simple and $L_K(E)$ is graded simple. By Theorem \ref{thm29}, graded ideals of $L_K(E)$ are generated by hereditary saturated sets in $E$, it follows that $E$ has no nonempty proper hereditary saturated set. 

Now, in contrary, suppose  $L_K(E)$ is not simple. By Theorem \ref{Alahmedi-Theorem2}, there exists simple subgraph $W$ such that every vertex outside $W$ is a balloon over $W$. Notice that if no such vertex exists, then $W=E$, a contradiction to the fact that $L_K(E)$ is not simple and $L_K(W)$ is simple by definition of simple subgraph. Hence, there exists $u\in E^0\setminus W$. By Theorem \ref{4thm3}, $E\setminus \{v\}$ is a hereditary saturated set in $E$. Thus, $E$ contains a proper hereditary saturated subset. Hence, $I(E\setminus \{v\})$ is a proper graded ideal of $L_K(E)$, a contradiction to $L_K(E)$ being graded simple. Consequently, $L_K(E)$ is simple.  \qed

Since simplicity directly implies graded simplicity, we obtain the following direct consequence of Theorem \ref{4thm1}. 
\begin{cor} \label{5cor1}
Let $E$ be a row-finite graph such that $[L_K(E), L_K(E)]$ is simple. Then $L_K(E)$ is graded simple if and only if $L_K(E)$ is simple.
\end{cor}

By Theorem \ref{thm29}, $L_K(E)$ being graded simple is equivalent to saying that $T_E$ is simple. Hence, together with  Corollary \ref{5cor1}, we have the following corollary.

\begin{cor}
Let $E$ be a row-finite graph such that $[L_K(E), L_K(E)]$ is simple. Then  $L_K(E)$ is simple if and only if $T_E$ is simple.
\end{cor}

\begin{thm} \label{4thm2}
Let $E$ be a finite  graph. Then the following is equivalent: 
\begin{enumerate}
    \item[{(i)}] $[L_K(E), L_K(E)]$ is simple and $L_K(E)$ is graded simple.
    \item[{(ii)}] $L_K(E)$ is simple and 
$1_{L_K(E)}= \sum_{v\in E^0} v\not \in [L_K(E), L_K(E)].$
\end{enumerate} 
\end{thm}
\pf Suppose $[L_K(E), L_K(E)]$ is simple and $L_K(E)$ is graded simple. Then by Theorem \ref{4thm1}, $L_K(E)$ is simple. By Theorem \ref{al1}, $$1_{L_K(E)}= \sum_{v\in E^0} v\not \in [L_K(E), L_K(E)].$$

Conversely, suppose $L_K(E)$ is simple and 
$$1_{L_K(E)}= \sum_{v\in E^0} v\not \in [L_K(E), L_K(E)].$$ Then by Theorem \ref{al1}, $[L_K(E), L_K(E)]$ is simple. Since simplicity implies graded simplicity. Thus, $L_K(E)$ is graded simple. \qed

The following corollary is also a direct consequence of Theorem \ref{thm29} applied on Theorem \ref{4thm2}.

\begin{cor}
Let $E$ be a finite  graph. Then the following is equivalent: 
\begin{enumerate}
    \item[{(i)}] $[L_K(E), L_K(E)]$ is simple and $T_E$ is simple.
    \item[{(ii)}] $L_K(E)$ is simple and 
$1_{L_K(E)}= \sum_{v\in E^0} v\not \in [L_K(E), L_K(E)].$
\end{enumerate}
\end{cor}

\section*{Acknowledgement} A.N. Sebandal wants to thank the Department of Science and Technology for a Ph.D. scholarship without which this paper would not be able to finish. She also wants to thank Technische Universit\"at Kaiserslautern of Germany for a research stay wherein this paper was finished.

\section*{Conflict of Interest} None of the authors has any conflict of interest in the conceptualization or publication of this work.


\end{document}